\documentclass [12pt]{article}
\usepackage{amsmath,amssymb}
\usepackage{amssymb,verbatim}

\usepackage{enumerate}

\newcommand{\interior}[1]{%
  {\kern0pt#1}^{\mathrm{o}}%
}
\usepackage{scalerel}
\usepackage{relsize}
\usepackage[T1]{fontenc}
\usepackage[tableposition=top]{caption}
\usepackage[countmax]{subfloat}

\usepackage{caption} 
\usepackage[nospace,noadjust]{cite}
\usepackage{float}
\newcommand\preceqdot{\mathrel{\ooalign{$\preceq$\cr
  \hidewidth\raise0.125ex\hbox{$\cdot\mkern0.5mu$}\cr}}}
\newcommand\precdott{\mathrel{\ooalign{$\prec$\cr
  \hidewidth\raise0.015ex\hbox{$\cdot\mkern0.5mu$}\cr}}} 
\usepackage{mathabx}
\usepackage{float}
\usepackage{blindtext}
\usepackage{titlesec}
\title{Sections and Chapters}
\usepackage{enumitem}
\usepackage{multirow}
\usepackage{array}
\newcolumntype{L}[1]{>{\raggedright\let\newline\\\arraybackslash\hspace{0pt}}m{#1}}
\newcolumntype{C}[1]{>{\centering\let\newline\\\arraybackslash\hspace{0pt}}m{#1}}
\newcolumntype{R}[1]{>{\raggedleft\let\newline\\\arraybackslash\hspace{0pt}}m{#1}}
\usepackage{blkarray}
\usepackage{amsmath}
\usepackage{tikz-cd} 
\usepackage{subfig}
\usepackage{faktor}
\usepackage{xfrac}
\usepackage{faktor} 
\usepackage{amsmath} 
\usepackage{amssymb}
\usepackage{mathtools}

\usepackage[utf8]{inputenc}
\usepackage{graphicx}
\usepackage{graphicx,amssymb,amstext,amsmath}
\usepackage{tikz}
\usepackage{amsmath}
\usepackage{amssymb}
\usepackage{amsthm}
\usepackage{hebcal}
\usetikzlibrary{arrows, decorations.markings}
\usepackage{blindtext}
\usepackage{amsthm}
\usepackage[T1]{fontenc}
\usepackage[utf8]{inputenc}
\usepackage{blindtext}
\usepackage[T1]{fontenc}
\usepackage[utf8]{inputenc}
\newcommand\Item[1][]{%
  \ifx\relax#1\relax  \item \else \item[#1] \fi
  \abovedisplayskip=0pt\abovedisplayshortskip=0pt~\vspace*{-\baselineskip}}
\usepackage{tikz}
\usetikzlibrary{positioning, arrows.meta}

\usepackage{cite}
\newtheorem{thm}{Theorem}[section]

\newtheorem{propx}{Proposition}

\newtheorem{thmx}{Theorem}
\newtheorem{innercustomgeneric}{\customgenericname}
\providecommand{\customgenericname}{}
\newcommand{\newcustomtheorem}[2]{%
  \newenvironment{#1}[1]
  {%
   \renewcommand\customgenericname{#2}%
   \renewcommand\theinnercustomgeneric{##1}%
   \innercustomgeneric
  }
  {\endinnercustomgeneric}
}

\newcustomtheorem{customthm}{Theorem}
\newcustomtheorem{customlemma}{Lemma}

\newtheorem{remark}[thm]{Remark}
\newtheorem{question}[thm]{Question}

\newtheorem{question
}{Question}
\newtheorem*{acknowledgement}{Acknowledgement}
\newenvironment{Proof}[1]{\par\noindent{\bf
    Proof{#1}:}\quad}{}
\def\0{{\bf 0}}
\def\a{{\bf a}}
\def\b{{\bf b}}
\def\d#1{{\rm d}#1\ }

\def\N{{\bf N}}
\def\P{{\bf P}}
\def\Q{{\bf Q}}
\def\R{{\bf R}}
\def\Z{{\bf Z}}

\def\B{{\rm B}}
\def\P{{\rm P}}

\def\MB{{\rm MB}}
\def\r{{\bf r}}
\def\f{{\rm f}}
\def\d{{\underline d}}
\def\HH{{\rm H}}
\def\T{{\bf T}}
\def\p{{\bf p}}

\def\longto{\mathop{\longrightarrow}\limits}

\usepackage{subfig}
\makeatletter
\def\keywords{\xdef\@thefnmark{}\@footnotetext}
\title {Exceptional Points for Density Modulo 1}
\date{}

\author{D. Berend\footnote{
Departments of Mathematics and Computer Science, Ben-Gurion
University, Beer Sheva 84105, Israel.
E-mail: berend@math.bgu.ac.il}
\footnote{Research supported in part by
the Milken Families Foundation Chair in
Mathematics.}
\and
M. D. Boshernitzan \footnote{Deceased; formerly of Department of Mathematics, Rice University, Houston, TX 77251, USA.}
\and G. Kolesnik\footnote{Department of Mathematics, California State University, Los Angeles, CA 90032, USA. E-mail: gkolesn@calstatela.edu}
\and
R. Kumar\footnote{Department of Mathematics, Ben-Gurion
University, Beer Sheva 84105, Israel. E-mail: kumarr@post.bgu.ac.il}
}
\begin{document}
\maketitle
\keywords{2020 {Mathematics Subject Classification:} Primary 11B05; 11J71; 11K06; 11K55.}
\keywords{{Key words and phrases. Distribution modulo 1, density modulo 1, Hausdorff dimension, box dimension, modified box dimension, IP-set.}}

\begin{abstract}
It is well known that almost every dilation of a sequence
of real numbers, that diverges to $\infty$, is dense modulo~1.
This paper studies the exceptional set of points -- those for which the dilation is not dense. Specifically, we consider the Hausdorff and modified box dimensions of the set of exceptional points. In particular, we show that the dimension of this set may be any number between 0 and 1. Similar results are obtained for two ``natural'' subsets of the set of exceptional points. Furthermore, the paper calculates the dimension of several sets of points, defined by certain constraints on their binary expansion.

\end{abstract}

\section{Introduction}
\medskip
Throughout this paper, $\r=(r_n)_{n=1}^\infty$ will denote a sequence of positive real numbers with $r_n\xrightarrow[n\to\infty]{} \infty$. It is well known that, if
\begin{equation} \label{condition}
r_{n+1}-r_n>C, \qquad n\geq 1,
\end{equation}
where $C>0$ is a constant, then the sequence $(r_n x)_{n=1}^\infty$ is
uniformly distributed modulo~1 for almost every real number~$x$ with respect to the Lebesgue measure (\cite{Weyl}; see also \cite [Ch.1, Corollary 4.3]
{kuipers-niederreiter}). However, the set of exceptional points may be quite large. In fact, denote
$$E_{\rm ud}(\r)=\{x\in \R\, : \, (r_nx)_{n=1}^\infty\, \textup{is not uniformly distributed modulo~1}\}.$$
Erd\H os and Taylor \cite[Theorem 14]{erdos-taylor} proved that, if $\r$ is an (increasing) integer-valued lacunary sequence, then $\dim_{\HH}E_{\rm ud}(\r)=1$, where $\dim_{\HH} A$ denotes the Hausdorff dimension of a set $A\subseteq\R$. Recall that $(r_n)_{n=1}^\infty$ is {\it lacunary} if $r_{n+1}/r_n\geq \lambda>1$ for each $n$. 

On the other hand, let $\r$ be a sequence of distinct integers, satisfying $|r_k|= O(k^{\beta +\varepsilon})$ for some fixed $\beta\geq 1$ and every $\varepsilon>0$. Then 
\begin{equation*}\label{result sec, dimension for ud}
\dim_{\HH}E_{\rm ud}(\r)\leq 1-\frac{1}{\beta}    
\end{equation*}
(see \cite{Sapiro, Salem, erdos-taylor}). Baker~\cite{Baker1} extended this result to real sequences. Ruzsa \cite[Theorem 6]{Ruzsa} proved, in the opposite direction, that for every $\beta \geq 1$ there exists an increasing sequence~$\r$ of integers, satisfying $r_k = O(k^{\beta})$, such that
\begin{equation*}\label{result sec, equality in dimension for ud}
\dim_{\HH}E_{\rm ud}(\r)= 1-\frac{1}{\beta}.    
\end{equation*}
For multi-dimensional analogs of $E_{\rm ud}(\r)$, we refer to \cite[Theorems 1,2]{Baker2}.

Here, we are interested in the set of exceptional points for density modulo 1. Denote
$$E(\r)=\{x\in \R\, : \, (r_nx)_{n=1}^\infty\, \textup{is not dense modulo~1}\}.$$
Note that $E(\r)\subseteq E_{\rm ud}(\r)$. Pollington and de Mathan, improving the above-mentioned result of Erd\H os and Taylor \cite[Theorem 14]{erdos-taylor}, relating to the exceptional set for uniform distribution, proved

\begin{thmx}\emph{\cite{pollington, de-mathan}}\label{theorem on lacunary}
If $(r_n)_{n=1}^\infty$ is lacunary, then $\dim_{\HH} E(\r)=1$.
\end{thmx}

Obviously, as $(r_n)_{n=1}^\infty$ increases more and more slowly, one would expect the set $E(\r)$ to become smaller and smaller. In view of Theorem \ref{theorem on lacunary}, it is natural to consider {\it sublacunary} sequences, namely sequences $(r_n)_{n=1}^\infty$ satisfying $r_{n+1}/r_n\longto_{n\to\infty}
1$. 

\begin{thmx}\emph{\cite{bosh2}}\label{theorem on sub-lacunary}
If $(r_n)_{n=1}^\infty$ is sublacunary, then $\dim_{\HH} E(\r)=0$.
\end{thmx}

Let $\r$ be the increasing sequence of integers consisting of all elements of a multiplicative subsemigroup $S$ of $(\N,\cdot)$. It is easy to see that, if $S$ is generated by a single integer, namely $S=\{q^n\, :\, n\geq 1\}$ for some integer $q\geq 2$, then $\dim_{\HH} E(\r)=1$. Moreover, Schmidt \cite{Schmidt} showed that
$$\dim_{\HH} \bigcap_{q\geq 2} E\left((q^n)_{n=1}^\infty\right)=1.$$

On the other hand, Furstenberg showed that, if $S$ contains two multiplicatively independent integers $q_1,q_2$ (i.e., $\log q_2/\log q_1$ is irrational), then  $E(\r)=\Q$ (\cite{Fustenberg2}; see also \cite[Theorem 1.1]{Bos5}).

We refer to \cite{bos3,bos2} for other relevant results on $E(\r)$ and to \cite{baker} for multi-dimensional analogs of $E(\r)$.

Theorems \ref{theorem on lacunary} and \ref{theorem on sub-lacunary}, taken together, hint at some kind of a 0-1 law, and thus raise

\begin{question}\emph{\cite{Peres}} \label{quetion1}\emph{Is $E(\r)$ always either of Hausdorff dimension 0 or of Hausdorff dimension 1?}
\end{question}

In Section \ref{main results}, we answer Question \ref{quetion1} negatively and also discuss the Hausdorff, box, and modified box dimensions of other sets related to $E(\r)$. We also obtain some results on the Hausdorff, box, and modified box dimension of certain sets of numbers defined by means of constraints on their binary expansion. In Section~\ref{Section Hausdoff dimension} we recall the basic properties of the various notions of dimension we use. Sections \ref{proof section of B}-\ref{Section on IP sets} are devoted to the proofs of the results.

\begin{acknowledgement}\emph{
We would like to thank B. Wang and J. Wu for helpful information, and in particular for their suggestion to consider the packing dimension also.}
\end{acknowledgement}

\section{Main Results}\label{main results}
\medskip

\subsection{The Dimension of Sets of Exceptional Points}\label{subsection 2.a}
\medskip
Our first result answers Question \ref{quetion1} negatively in a strong way.

\begin{thm} \label{E-may-be-any-dim}
For every $0\leq d\leq 1$, there exists a sequence of integers
$(r_n)_{n=1}^\infty$ such that 
$$\dim_{\HH} E(\r)= {\underline{\dim}}_{\MB} E(\r)= d,$$
where ${\underline{\dim}}_{\MB}$ denotes modified lower box dimension.
\end{thm}

\begin{remark}\emph{Since $\r$ consists of integers,
the set $E(\r)$ contains all rationals, so trivially ${\underline{\dim}}_{\B} E(\r)=1$, where ${\underline{\dim}}_{\B}$ denotes lower box dimension.
}
\end{remark}

The set $E(\r)$ contains a subset $E_0(\r)$, consisting of all those $x\in\R$ for which $r_nx \xrightarrow[n\to\infty]{} 0\; (\textup{mod}\, 1)$. Obviously,
$E_0(\r)$ is an additive subgroup of $\R$. It follows in
particular from \cite [Theorem 7.1] {bosh2} that, if the
sequence of ratios $(r_{n+1}/r_n)_{n=1}^\infty$ is bounded above,
then $E_0(\r)$ is countable. (However, Dubickas \cite{Dubickas} and Bugeaud \cite{Bugeaud} showed that there exist ``slowly increasing'' sequences $\r$ of integers for which $E(\r)$ contains irrationals.) On the other hand, if $\r$ is integer-valued and $r_{n+1}/r_n\xrightarrow[n\to\infty]{} \infty$, then
$\dim_\HH E_0(\r)=1$ (see \cite[Theorem 8]{erdos-taylor}). Similarly to Theorem
\ref{E-may-be-any-dim}, for general sequences~$\r$,
this dichotomy is not true any more. In fact, let $\alpha$ be any irrational and $(q_n)_{n=1}^\infty=(q_n(\alpha))_{n=1}^\infty$ be the sequence of the denominators of the convergents in the continued fraction expansion of $\alpha$. Wang and Wu \cite[Corollary 3.3]{Wang1} proved that for each $0\leq d\leq 1$, there exists an irrational $\alpha \in [0,1)$ such that $\dim_{\HH}E_{0}\left((q_n)_{n=1}^\infty\right)=d$.

The group $E_0(\r)$ in turn has a subgroup $E_\f(\r)$,
consisting of all points $x\in [0,1)$ for which
$\sum_{n=1}^\infty \|r_n x\|<\infty$, where $\|t\|$ denotes the
distance of a real number $t$ from the nearest integer. This set played a central role in Erd\H os-Taylor's paper
\cite{erdos-taylor}. It follows from \cite[Theorem 6]{erdos-taylor} that, for each $0\leq d\leq 1$, there exits an integer-valued sequence~$\r$ such that $\dim_{\HH}E_{\f}(\r)=d$. Wang and Wu \cite[Theorem 1.2]{Wang2} proved a strengthened version of that result. Thus, it follows from \cite{erdos-taylor, Wang1, Wang2} that, for every $0\leq d\leq 1$, one can find integer-valued sequences $\r$ and $\r'$ such that $E_{0}(\r)=d$ and $E_{f}(\r')=d$. The following result shows that one may have the same sequence $\r$ for both equalities. Moreover, the proof provides a very simple such sequence $\r$, having transparent exceptional sets $E_0(\r)$ and $E_\f(\r)$ of the same Hausdorff dimension.

\begin{thm} \label{E_0-may-be-any-dim}
For every $0\leq d\leq 1$, there exists a sequence of integers
$(r_n)_{n=1}^\infty$ such that 
\begin{description}
\item {1.} \vspace{-0.3 cm} $$\dim_{\HH} E_0(\r)= {\underline{\dim}}_{\MB} E_0(\r)=d.$$
\item {2.} $$\dim_{\HH} E_\f(\r)= {\underline{\dim}}_{\MB} E_\f(\r)=d.$$
\end{description}
\end{thm}

\begin{remark}\emph{
The groups $E_{0}(\r)$ and $E_\f(\r)$, for the sequence $\r$ we construct, contain the group of all dyadic rationals in $\T$, and therefore ${\underline{\dim}}_{\B} E_0(\r)= {\underline{\dim}}_{\B} E_\f(\r)=1$.
}
\end{remark}
\vspace{1.cm}

\subsection {The Dimension of Sets Defined by Binary~Expansions}\label{subsection 2.b}
\medskip
The following theorems will be used to prove the results in Subsection \ref{subsection 2.a}, but they are of independent interest. (Some of them may be known for some notions of dimension, but to be self-contained, we present their proofs.)
Recall that the {\it lower} and {\it upper} {\it density} of a set $S\subseteq\N$, respectively, are
given by
\begin{equation}
\d(S)=\liminf_{N\to\infty} \frac{\left|S\cap
[1,N]\right|}{N},
\label{lower density}
\end{equation}
\begin{equation}
\bar{d}(S)=\limsup_{N\to\infty} \frac{\left|S\cap
[1,N]\right|}{N}.
\label{upper density}
\end{equation}

\begin{thm} \label{hausdorff-dim-calc}
Let $S\subseteq\N$, and let $X(S)\subseteq [0,1]$ be the set of
all numbers $x$, in whose infinite binary expansion $x=\sum_{n=1}^\infty\xi_n 2^{-n}$ we have $\xi_n=0$ for each $n\notin
S$. Then:
\begin{description}
\item {1.} \vspace{-0.4 cm} $$\dim_{\HH} X(S)={\underline{\dim}}_{\MB} X(S)={\underline{\dim}}_{\B} X(S)=\d(S).$$
\item {2.}  \vspace{-0.4 cm} $$\dim_{\P} X(S)={\overline{\dim}}_{\MB} X(S)={\overline{\dim}}_{\B} X(S)= \bar{d}(S),$$
\end{description}
where ${\overline{\dim}}_{\B} F$, ${\overline{\dim}}_{\MB} F$, and $\dim_{\P}F$ are the upper box, modified upper box, and packing dimensions of $F$, respectively.
\end{thm}

\begin{remark}\label{remark on X(a,b)}\emph{In other words, $X(S)$ is the set of numbers whose infinite binary expansion is ``free'' at $S$ and ``restricted'' outside $S$; digits at places belonging to $S$ are allowed to be either $0$ or $1$, while the other digits must be $0$. The theorem has an intuitive appeal: The size of $X(S)$, measured by (various versions of) dimension, is the size of $S$, measured by (the corresponding versions of) density.}
\end{remark}

As the theorem is trivial if either $S$ or its complement
in $\N$ is finite, we shall implicitly assume that both sets are infinite. The sets we shall look at will typically consist of large blocks of consecutive integers, with large holes between these blocks. Thus, let us reformulate the theorem in a slightly different way.

\begin{thm}\label{hausdorff-dim-calc'}
Let $\a=(a_i)_{i=1}^\infty$ and $\b=(b_i)_{i=1}^\infty$ be sequences of integers such that 
$$0=a_1\leq b_1<a_2<b_2<\cdots<a_i<b_i<\cdots.$$
Let $X(\a,\b)$ be the subset of $[0,1]$, consisting of all those
numbers $x$, whose infinite binary expansion
$x=\sum_{n=1}^\infty\xi_n 2^{-n}$ satisfies the condition
\begin{equation*} \label{constraint}
\xi_n=0, \qquad n\in \bigcup_{i=1}^\infty [b_i+1,a_{i+1}].
\end{equation*}
Then:
\begin{description}
\item {1.} \vspace{-0.6 cm} \begin{equation} \label{1st-formula-for-dim}
\dim_{\HH} X(\a,\b)={\underline{\dim}}_{\MB} X(\a,\b)={\underline{\dim}}_{\B} X(\a,\b)=
    \liminf_{k\to\infty} \frac{\sum_{i=1}^k (b_i-a_i)}{a_{k+1}}.
\end{equation}
\item {2.} \vspace{-0.6 cm}
\begin{equation}\label{result sec, formula for packing dim of X i}
\dim_{\P}X(\a,\b)={\overline{\dim}}_{\MB} X(\a,\b)= {\overline{\dim}}_{\B} X(\a,\b)= \limsup_{k\to\infty} \frac{\sum_{i=1}^k (b_i-a_i)}{b_{k}}.    
\end{equation}
\end{description}
\end{thm}

\begin{thm} \label{hausdorff-dim-calc''}
Let $\a=(a_i)_{i=1}^\infty$ and $\b=(b_i)_{i=1}^\infty$ be as
in Theorem \ref{hausdorff-dim-calc'}.
Let $X'(\a,\b)$ be the subset of $[0,1]$, consisting of all those numbers $x$, whose infinite binary expansion
$x=\sum_{n=1}^\infty\xi_n 2^{-n}$ satisfies the condition
\begin{equation*} \label{constraint'}
\xi_{b_i+1}=\xi_{b_i+2}=\cdots=\xi_{a_{i+1}}, \qquad
i=1,2,3,\ldots.
\end{equation*}
Then:
\begin{description}
\item {1.} \vspace{-0.6 cm}\begin{equation}\label{proposition on X' equation}
\dim_{\HH} X'(\a,\b)={\underline{\dim}}_{\MB} X'(\a,\b)={\underline{\dim}}_{\B} X'(\a,\b)=
    \liminf_{k\to\infty} \frac{\sum_{i=1}^k
(b_i-a_i+1)}{a_{k+1}}.    
\end{equation}
\item {2.} \vspace{-0.8 cm}  
\begin{equation}\label{proposition on X' equation of packing dim}
\dim_{\P}X'(\a,\b)={\overline{\dim}}_{\MB} X'(\a,\b)={\overline{\dim}}_{\B} X'(\a,\b)= \limsup_{k\to\infty} \frac{\sum_{i=1}^k (b_i-a_i+1)}{b_{k}}.  
\end{equation}
\end{description}

\end{thm}

\begin{remark}\label{remark on X'(a,b)}\emph{
Similarly to Remark \ref{remark on X(a,b)}, here $X'(\a,\b)$ is the set of numbers whose infinite binary expansion is ``free'' at $\bigcup_{i=1}^\infty [a_{i}+1,b_{i}]$ and ``restricted'' in $\bigcup_{i=1}^\infty [b_{i}+1,a_{i+1}]$; digits at places belonging to $\bigcup_{i=1}^\infty [a_{i}+1,b_{i}]$ are allowed to be either $0$ or $1$, while the digits in each block $[b_{i}+1,a_{i+1}]$ for $i\geq 1$, are either all $0$ or all $1$. The numerator $\sum_{i=1}^k (b_i-a_i)$ on the right-hand side of \eqref{1st-formula-for-dim} and \eqref{result sec, formula for packing dim of X i} counts the ``degrees of freedom'' up to the $b_k$-th place when choosing a number in $X(\a,\b)$. For $X'(\a,\b)$, each block $[b_{i}+1,a_{i+1}]$ provides an additional degree of freedom. This explains the additional +1s on the right-hand side of \eqref{proposition on X' equation} and \eqref{proposition on X' equation of packing dim} on the top of the right-hand side of  \eqref{1st-formula-for-dim} and \eqref{result sec, formula for packing dim of X i}.
}
\end{remark}

\section{Various Notions of Dimension}\label{Section Hausdoff dimension}
\medskip
In this section, we recall the definitions and basic properties of Hausdorff, box, and modified box, dimensions. We also mention briefly the packing dimension.
\subsection{Hausdorff Dimension}
Let $U$ be a non-empty subset of $\R^n$. The {\it diameter} of $U$ is defined as
$$|U|= \sup\{|x-y|\, :\, x,y\in U\}.$$
Let $\delta>0$. A countable collection $\{U_i\}_{i=1}^\infty$ of subsets of $\R^n$ is a $\delta$-{\it cover} of $F\subseteq \R^n$ if 
$$F\subseteq \bigcup_{i=1}^\infty U_i,\qquad  0\leq |U_i|\leq \delta,\, i=1,2,\ldots.$$
Let $s$ be a non-negative real number. Denote:
$$\mathcal{H}_{\delta}^s(F)=\inf\left\{\sum_{i=1}^\infty|U_i|^s\, :\, \{U_i\}_{i=1}^\infty,  \textup{is a}\, \delta\textup{-cover of}\, F\right\},\qquad \delta>0.$$
As $\delta$ decreases, $\mathcal{H}_{\delta}^s(F)$ increases. The $s$-{\it dimensional Hausdorff measure} of $F$ is given by
$$\mathcal{H}^s(F)= \lim_{\delta\to 0}\mathcal{H}_{\delta}^s(F),$$
and the {\it Hausdorff dimension} of $F$ by
$$\dim_{\HH} F= \inf\{s\, :\, \mathcal{H}^s(F)=0\}=\sup\{s\, :\, \mathcal{H}^s(F)=\infty\}.$$
The Hausdorff dimension satisfies the following properties (see, for example, \cite[Proposition 4.3]{Carlo}):
\begin{enumerate}
\item {\it Monotonicity.} If $E\subseteq F$, then $\dim_{\HH} E\leq \dim_{\HH} F$.
\item {\it Countable stability.} If $F_1, F_2,\ldots$ is a sequence of sets, then
\begin{equation}\label{Hausdorff dim of union}
\dim_{\HH} \bigcup_{i=1}^\infty F_i= \sup\{\dim_{\HH} F_i\, :\, i=1,2,\ldots\}.    
\end{equation}
\item {\it Lipschitz invariance.} If $f$ is a bi-Lipschitz transformation, then $\dim_{\HH} f(F)=\dim_{\HH} F$.
\end{enumerate}
\subsection{Box Dimension}
Let $F$ be a non-empty bounded subset of $\R^n$. For $\varepsilon>0$, let $N_{\varepsilon}(F)$ be the smallest number of sets of diameter at most $\varepsilon$ which can cover $F$. The {\it lower box dimension} and the {\it upper box dimension} of~$F$ are defined by
\begin{equation*}
  {\underline{\dim}}_{\B} F   = \liminf_{\varepsilon\to 0}\frac{\log N_{\varepsilon}(F)}{-\log \varepsilon},
  \end{equation*}    
\begin{equation*}
   {\overline{\dim}}_{\B} F = \limsup_{\varepsilon\to 0}\frac{\log N_{\varepsilon}(F)}{-\log \varepsilon},
\end{equation*}   
respectively.
If the two are equal, we refer to the common value as the {\it box dimension} of~$F$ and denote it by $\dim_{\B} F$. The lower and the upper box dimensions are monotonic and Lipschitz invariant.
It follows from the definitions of Hausdorff and lower box dimensions that 
$$\dim_{\HH} F \leq {\underline{\dim}}_{\B}F$$
(see, for example \cite[(3.17)]{falconer}).

The following proposition provides upper bounds on the lower and the upper box dimensions. 
\begin{propx}{\emph{\cite[Proposition 4.1]{falconer}}}\label{Proposition 4.1 falconer}
Suppose $F$ can be covered by $n_k$ sets of diameter at most $\varepsilon_k$ with $\varepsilon_k\to 0$ as $k\to \infty$. Then:
$$\dim_{\HH} F \leq {\underline{\dim}}_{\B}F \leq \liminf_{k\to \infty} \frac{\log n_k}{-\log \varepsilon_k}.$$
If, moreover, $\varepsilon_{k+1}>c \varepsilon_k$ for some $c>0$, then
$${\overline{\dim}}_{\B} F \leq \limsup_{k\to \infty} \frac{\log n_k}{-\log \varepsilon_k}.$$
\end{propx}
A {\it mass distribution} on $F$ is a measure whose support is contained in $F$, such that $0< \mu(F)<\infty$.

\begin{propx}{\emph{\cite[Mass distribution principle 4.2]{falconer}}}\label{Proposition mass distribution falconer}
Let $\mu$ be a mass distribution on $F$ and suppose that for some $s$ there are numbers $c>0$ and $\varepsilon>0$ such that
$$\mu(U)\leq c|U|^s,\qquad |U|\leq \varepsilon.$$
Then: 
$$s \leq \dim_{\HH} F \leq {\underline{\dim}}_{\B}F\leq \overline{\dim_{\B}}F.$$
\end{propx}

The upper box dimension is finitely stable 
$${\overline{\dim}}_{\B} (E\cup F) = \max\, \{{\overline{\dim}}_{\B} E, {\overline{\dim}}_{\B} F\},\qquad E,F \subseteq \R^n,$$
but the lower box dimension is not \cite[p.48]{falconer}. The lower and upper box dimension are invariant under closure \cite[Proposition 3.4]{falconer}. Note that countable closed sets may have positive box dimension. For example, if $F=\left\{0,1,\frac{1}{2},\frac{1}{3},\ldots\right\}$, then $\dim_{\B} F=1/2$ (see \cite[Example 3.5]{falconer}).

\subsection{Modified Box Dimension}

To overcome the peculiarities of box dimension, a notion of modified box dimension has been introduced. Let $F\subseteq \R^n$. The modified lower and upper box dimensions of $F$ are defined by
\begin{equation*}
{\underline{\dim}}_{\MB}F= \inf \left\{\sup_{i} {\underline{\dim}}_{\B}F_i\, :\, F\subseteq \bigcup_{i=1}^\infty F_i \right\},
\end{equation*}

\begin{equation*}
{\overline{\dim}}_{\MB}F= \inf \left\{\sup_{i} {\overline{\dim}}_{\B}F_i\, :\, F\subseteq \bigcup_{i=1}^\infty F_i \right\},
\end{equation*}
respectively, where in both cases the infimum is over all countable covers $\{F_i\}_{i=1}^\infty$ of $F$. Obviously, for every $F\subseteq \R^n$,
\begin{equation}\label{pre, relation in box and modified box}
{\underline{\dim}}_{\MB}F\leq {\underline{\dim}}_{\B}F,\qquad {\overline{\dim}}_{\MB}F\leq {\overline{\dim}}_{\B}F,    
\end{equation}
and
\begin{equation}\label{pre, relation in Hausdoff and modified}
 \dim_{\HH} F\leq {\underline{\dim}}_{\MB}F  \leq {\overline{\dim}}_{\MB}F \leq {\overline{\dim}}_{\B}F. 
\end{equation}
The modified upper and lower box dimensions are monotonic, countably stable, and Lipschitz invariant. In the case of compact set, the following proposition gives a sufficient condition for the box and the  modified box dimension to coincide.

\begin{propx}{\emph{\cite[Proposition 3.6]{falconer}}}\label{Proposition 3.6 falconer}
Let $F\subseteq \R^n$ be compact. Suppose that
$${\overline{\dim}}_{\B}(F\cap V)= {\overline{\dim}}_{\B}F$$
for all open sets $V$ that intersect $F$. Then ${\overline{\dim}}_{\MB} F={\overline{\dim}}_{\B} F$.
\end{propx}
\vspace{1 cm}

Note that the Hausdorff dimension is defined using measures, whereas neither the box dimension nor the modified box dimension is based on measures. Another notion of dimension defined via measures is the {\it packing dimension}, denoted by $\dim_{\P}$ (see \cite[Section 3.4]{falconer}). As the packing dimension coincides with the modified upper box dimension for $F\subseteq \R^n$ \cite[Prposition 3.8]{falconer}, we do not present the definition of the former.

We refer to \cite{Carlo, falconer} for more details on Hausdorff, box, modified box, and packing dimensions.

\section{Proofs of Theorems \ref{hausdorff-dim-calc'} and \ref{hausdorff-dim-calc''}}\label{proof section of B}
\medskip
\begin{Proof}{\ of Theorem \ref{hausdorff-dim-calc'}}
1. Denote
$$d_1= \liminf_{k\to\infty} \frac{\sum_{i=1}^k (b_i-a_i)}{a_{k+1}}.$$
Let $D_k$ denote the following finite set of dyadic rationals:
\begin{equation*}
 \begin{split}
D_k= \left\{x= \sum_{n=1}^{b_k}\frac{\xi_n}{2^n}\, :\, \xi_n\in \{0,1\}\, \forall\,  n,\, \xi_n=0\,  \forall\,  n\in \bigcup_{i=1}^{k-1}[b_{i}+1, a_{i+1}]\right\},\qquad k\geq 1.  
 \end{split}   
\end{equation*}
Clearly, 
$$D_k\subseteq X(\a,\b)\subseteq \bigcup_{x\in D_k}
[x,x+2^{-a_{k+1}}),\qquad k\geq 1.$$
The cardinality of $D_k$ is $2^{\sum_{i=1}^k(b_i-a_i)}$. Thus, for each $k$, the set $X(\a,\b)$ is contained in a union of $2^{\sum_{i=1}^k (b_i-a_i)}$ intervals of length $2^{-a_{k+1}}$ each, which yields, by Proposition \ref{Proposition 4.1 falconer},
\begin{equation} \label{upper-bound}
\dim_{\HH} X(\a,\b)\leq {\underline{\dim}}_{\B} X(\a,\b)\le 
    \liminf_{k\to\infty} \frac{\log 2^{\sum_{i=1}^k (b_i-a_i)}}
                              {-\log 2^{-a_{k+1}}}
      =\liminf_{k\to\infty} \frac{\sum_{i=1}^k
(b_i-a_i)}{a_{k+1}}=d_1.
\end{equation}

In the other direction, we first define a probability measure
$\mu$ on $[0,1)$ as follows. For each $n\in\N$, let

\begin{equation*}
\mu_n=\left\{
          \begin{array}{ll}
            \delta_0, \qquad & n\in \bigcup_{i=1}^\infty [b_i+1,a_{i+1}],\\
            \frac{1}{2}\delta_0+\frac{1}{2}\delta_{2^{-n}},
              \qquad & n\in \bigcup_{i=1}^\infty [a_i+1,b_{i+1}],
          \end{array} \right.    
\end{equation*}
where $\delta_t$ denotes a point mass at $t$.
Now let $\mu=*_{n=1}^\infty \mu_n$ be the convolution of all
$\mu_n$-s. Clearly, $\mu$ is supported on $X(\a,\b)$.

We would like to employ Proposition \ref{Proposition mass distribution falconer}. This requires upbounding the $\mu$-measure of sets of a small diameter. It clearly suffices to deal with intervals. Thus, let $I$ be a small interval. Assume first that $I$ is of the form $[l/2^n,(l+1)/2^n]$ for some $n$ and $l\leq 2^n-1$. We have
\begin{itemize}
\item If $l/2^n \notin \bigcup_{i=1}^\infty D_i$, then 
\begin{equation}\label{proof sec, equations1 for measure I for X}
 \mu(I)=0.   
\end{equation}
\item If $l/2^n \in \bigcup_{i=1}^\infty D_i$ and $a_{k}+1\leq n\leq b_{k}$ for some $k$, then 
\begin{equation}\label{proof sec, equations3 for measure I for X}
\mu(I)= \displaystyle \frac{1}{2^{\sum_{i=1}^{k-1}(b_i-a_i) + (n-a_{k})}}.  
\end{equation}
\item If $l/2^n \in \bigcup_{i=1}^\infty D_i$ and $b_{k}+1\leq n\leq a_{k+1}$ for some $k$, then 
\begin{equation}\label{proof sec, equations2 for measure I for X}
\mu(I)= \displaystyle \frac{1}{2^{\sum_{i=1}^{k}(b_i-a_i)}}.  
\end{equation}
\end{itemize}
Let $\varepsilon>0$. The conditions of the theorem ensure that
$$\sum_{i=1}^{k-1} (b_i-a_i)\geq (d_1-\varepsilon) a_{k},$$
for all sufficiently large $k$.
Hence, for $n\geq a_{k}+1$,
\begin{equation}\label{proof sec, equation for inequality of sum in terms of eps}
\begin{split}   
\sum_{i=1}^{k-1} (b_i-a_i) + (n-a_{k})&\geq n- (1- d_1+\varepsilon)a_{k} \geq n- (1- d_1+\varepsilon)n= (d_1-\varepsilon)n.\\
\end{split}    
\end{equation}
By \eqref{proof sec, equations1 for measure I for X}-\eqref{proof sec, equation for inequality of sum in terms of eps}, we have
\begin{equation}\label{proof sec, equation of measure of I'}
\mu(I)\leq \displaystyle \frac{1}{2^{\sum_{i=1}^{k-1}(b_i-a_i) + (n-a_{k})}}\leq \frac{1}{2^{n(d_1-\varepsilon)}}=|I|^{d_1-\varepsilon}.    
\end{equation}
Now let $I$ be any interval with $\frac{1}{2^{m+1}}\leq |I|<\frac{1}{2^{m}}$ for a large $m$. For a suitable $0\leq l\leq 2^m-2$, we have
\begin{equation}\label{proof sec, equation I contained in I1 and I2}
I\subseteq \left[\frac{l}{2^m}, \frac{l+1}{2^m}\right]\bigcup \left[\frac{l+1}{2^m}, \frac{l+2}{2^m}\right]=I_1 \cup I_2.
\end{equation}
Using \eqref{proof sec, equation of measure of I'} and \eqref{proof sec, equation I contained in I1 and I2},
\begin{equation}\label{proof sec, equation for measure of I}
\mu(I)\leq \mu(I_1) + \mu(I_2)\leq \left(\frac{1}{2^m}\right)^{d_1-\varepsilon} + \left(\frac{1}{2^{m}}\right)^{d_1-\varepsilon}= 2^{1+d_1-\varepsilon} \left(\frac{1}{2^{m+1}}\right)^{d_1-\varepsilon} \leq 2^{1+d_1}|I|^{d_1-\varepsilon}.
\end{equation}
By \eqref{proof sec, equation for measure of I} and Proposition \ref{Proposition mass distribution falconer},
$$\dim_{\HH} X(\a,\b)\ge d_1-\varepsilon,$$
and therefore
\begin{equation} \label{lower-bound}
\dim_{\HH} X(\a,\b)\ge d_1.
\end{equation}
Using \eqref{pre, relation in box and modified box}, \eqref{pre, relation in Hausdoff and modified}, \eqref{upper-bound}, and \eqref{lower-bound}, we get
$$d_1\leq \dim_{\HH} X(\a,\b) \leq {\underline{\dim}}_{\MB} X(\a,\b)\leq {\underline{\dim}}_{\B}X(\a,\b) \leq d_1,$$
which completes the proof of this part.
\vspace{1 cm}

\noindent 2. 
Denote
$$d_2=  \limsup_{k\to\infty} \frac{\sum_{i=1}^k (b_i-a_i)}{b_{k}}.$$
We divide the proof into three parts:
\begin{description}
\item {$(i)$} \vspace{-0.3 cm} $${\overline{\dim}}_{\B} X(\a,\b)=d_2.$$
\item {$(ii)$}  \vspace{-0.3 cm} $${\overline{\dim}}_{\MB} X(\a,\b)=d_2.$$
\item {$(iii)$}  \vspace{-0.3 cm} $$\dim_{\P} X(\a,\b)=d_2.$$
\end{description}
\vspace{0.25 cm}

\noindent $(i)$ For any $n$, cover $X(\a,\b)$ by the union of all dyadic intervals $\left[\frac{l}{2^n}, \frac{l+1}{2^n}\right]$, $0\leq l\leq 2^n-1$, intersecting $X(\a,\b)$. These intervals correspond to numbers $l$ whose binary expansion has 0-digits at all places in $\bigcup_{i=1}^\infty[b_{i}+1,a_{i+1}]$. The number of these intervals is 
\begin{equation}\label{equation upper box dimension Nn}
\begin{split}
M_{2^{-n}}&= \begin{dcases}
   2^{(b_1-a_1)+ \cdots + (b_{k-1}-a_{k-1})+ (n- a_k)}, &\qquad a_{k}+ 1\leq n\leq b_k,\\
2^{(b_1-a_1)+ \cdots + (b_k- a_k)}, &\qquad b_{k}+ 1\leq n\leq a_{k+1}.\\
\end{dcases}
\end{split}   
\end{equation}
Therefore, 
\begin{equation}\label{equation upper box dimension Nn/en1}
\begin{split}
\frac{\log M_{2^{-n}}}{-\log 2^{-n}}
&= \begin{dcases}
    \frac{(b_1-a_1)+ \cdots + (b_{k-1}-a_{k-1})+ (n- a_k)}{n},  \qquad &a_{k}+ 1\leq n\leq b_k,\\
  \frac{(b_1-a_1)+ \cdots + (b_k- a_k)}{n}, \qquad&b_{k}+ 1\leq n\leq a_{k+1}.\\
\end{dcases}\\
\end{split}   
\end{equation}
It follows from \eqref{equation upper box dimension Nn/en1} that the sequence $\left(\frac{\log M_{2^{-n}}}{-\log 2^{-n}}\right)_{n=1}^\infty$ is increasing in the interval $[a_{k}, b_{k}]$ and decreasing in $[b_{k}, a_{k+1}]$. By Proposition~\ref{Proposition 4.1 falconer} and \eqref{equation upper box dimension Nn/en1}, we have
\begin{equation}\label{equation upper box dimension bound from above}
 {\overline{\dim}}_{\B} X(\a,\b) \leq \limsup_{n\to \infty} \frac{\log M_{2^{-n}}}{-\log 2^{-n}}=\limsup_{k\to \infty} \frac{\sum_{i=1}^k(b_i-a_i)}{b_k}=d_2.   
\end{equation}

In the other direction, we first note that every two elements of $D_k$ are at least $2^{-b_k}$ apart. Hence, each subinterval of $[0,1]$ of length $2^{-b_k}$ contains at most two points of $D_k$. Since $$|D_k|=2^{\sum_{i=1}^k(b_i-a_i)},$$
the minimal number of intervals of length $2^{-b_k}$ required to cover $D_k$ is at least $\frac{1}{2}\cdot 2^{\sum_{i=1}^k(b_i-a_i)}$.
Let $N_{2^{-b_k}}$ be the minimal number of intervals of length $2^{-b_k}$ covering the set $X(\a,\b)$. Since $D_k\subseteq X(\a,\b)$, we have 
\begin{equation*}
N_{2^{-b_k}}\geq \frac{1}{2}\cdot 2^{\sum_{i=1}^k(b_i-a_i)}.   
\end{equation*}
and therefore:
$$\frac{\log N_{2^{-b_k}}}{-\log 2^{-b_k}} \geq \frac{\log 2^{-1} 2^{\sum_{i=1}^k(b_i-a_i)}}{-\log 2^{-b_k}}.$$
Hence,
\begin{equation}\label{equation upper box dimension bound from below of Nn/en}
{\overline{\dim}}_{\B} X(\a,\b)= \limsup_{\delta \to 0} \frac{\log N_{\delta}}{-\log \delta }\geq \limsup_{k\to \infty} \frac{\log N_{2^{-b_k}}}{-\log 2^{-b_k}} \geq  \limsup_{k\to \infty} \frac{\sum_{i=1}^k(b_i-a_i)}{b_k}=d_2.
\end{equation}
By \eqref{equation upper box dimension bound from above} and
\eqref{equation upper box dimension bound from below of Nn/en},
\begin{equation}\label{proof sec, equation of upper box equal to d2 for X}
{\overline{\dim}}_{\B} X(\a,\b)= \limsup_{k\to \infty} \frac{\sum_{i=1}^k(b_i-a_i)}{b_k}=d_2,    
\end{equation}
which proves this part.
\vspace{0.25 cm}

\noindent $(ii)$
Let $V$ be any open set such that $X(\a,\b)\bigcap V\neq \varnothing$. By the monotonocity of ${\overline{\dim}}_{\B}$, we have
\begin{equation}\label{proof sec, upper box dim of XcapV lessthan X}
{\overline{\dim}}_{\B} X(\a,\b)\bigcap V\leq {\overline{\dim}}_{\B} X(\a,\b).    
\end{equation}
Let $x= \sum_{i=1}^\infty\frac{\xi_i}{2^i} \in X(\a,\b)\bigcap V$. Then for sufficiently large $k$,
$$\bar{x}= \sum_{i=1}^k \frac{\xi_i}{2^i}\in X(\a,\b)$$
and 
$$\left[\bar{x}, \bar{x} + 1/2^k\right] \subseteq V.$$
Hence 
\begin{equation}\label{proof sec, upper box of XcapV>X[x,x+1/2^k]}
 {\overline{\dim}}_{\B} X(\a,\b)\bigcap V \geq {\overline{\dim}}_{\B} X(\a,\b)\bigcap\left[\bar{x}, \bar{x} + 1/2^k\right].  
\end{equation}
Now all non-empty sets of the form
$$X(\a,\b)\bigcap \left[l/2^k,(l+1)/2^k\right], \qquad (0\leq l\leq 2^k-1),$$
are translates of each other. Since ${\overline{\dim}}_{\B}$ is finitely stable, we get
\begin{equation}\label{proof sec, upper box of X greater than XcapI}
\begin{split}
{\overline{\dim}}_{\B} X(\a,\b)&=\max\left\{{\overline{\dim}}_{\B} X(\a,\b) \bigcap \left[\frac{l}{2^k},\frac{l+1}{2^k}\right]\, :\, 0\leq l\leq 2^k-1\right\}\\
&={\overline{\dim}}_{\B} X(\a,\b)\bigcap\left[\bar{x}, \bar{x} + \frac{1}{2^k}\right].
\end{split}    
\end{equation}
By \eqref{proof sec, upper box dim of XcapV lessthan X}-\eqref{proof sec, upper box of X greater than XcapI},
\begin{equation}\label{proof sec, equaltion upper box of Xcap V equal to X}
{\overline{\dim}}_{\B} X(\a,\b)\bigcap V=  {\overline{\dim}}_{\B} X(\a,\b).
\end{equation}
Since $V$ is an arbitary open set intersecting $X(\a,\b)$, by Proposition \ref{Proposition 3.6 falconer}, \eqref{proof sec, equation of upper box equal to d2 for X}, and \eqref{proof sec, equaltion upper box of Xcap V equal to X},
\begin{equation}\label{proof sec, equation of modified upper of X equal to d2}
{\overline{\dim}}_{\MB} X(\a,\b)=  {\overline{\dim}}_{\B} X(\a,\b)=d_2,  \end{equation}
which completes the proof of this part.
\vspace{0.25 cm}

\noindent $(iii)$ By \cite[Proposition 3.8]{falconer}, 
\begin{equation}\label{proof sec. equation of dim(H)=dim(MB)}
\dim_{\P}F= {\overline{\dim}}_{\MB} F, \qquad F\subseteq \R^n.   
\end{equation}
By \eqref{proof sec, equation of modified upper of X equal to d2} and \eqref{proof sec. equation of dim(H)=dim(MB)},
\begin{equation}\label{proof sec, equation of modified upper box equal to packing X}
\dim_{\P}X(\a,\b)= {\overline{\dim}}_{\MB} X(\a,\b)=d_2,   
\end{equation}
which completes the proof of this part.
\end{Proof}
\vspace{1.0 cm}

\begin{Proof}{\ of Theorem \ref{hausdorff-dim-calc''}}
1. Denote 
$$d_1'=\liminf_{k\to\infty} \frac{\sum_{i=1}^k
(b_i-a_i+1)}{a_{k+1}}.$$
For the upper bound, we proceed in the same way as in Theorem \ref{hausdorff-dim-calc'}.1, counting the number of dyadic intervals of each length intersecting our set. In fact, consider the following finite sets of dyadic rationals:
\begin{equation*}
 \begin{split}
D_k'= \left\{x= \sum_{n=1}^{a_{k+1}}\frac{\xi_n}{2^n}\, :\, \xi_{b_i+1}=\xi_{b_i+2}=\cdots= \xi_{a_{i+1}},\, 1\leq i\leq k\right\},\qquad k\geq 1.  
 \end{split}   
\end{equation*}
Clearly, 
$$D_{k}' \subseteq X'(\a,\b)\subseteq \bigcup_{x\in D_k'}
[x,x+2^{-a_{k+1}}),\qquad k\geq 1.$$
The cardinality of $D_k'$ is $2^{\sum_{i=1}^k(b_i-a_i+1)}$. Thus, for each $k$, the set $X'(\a,\b)$ is
contained in a union of $2^{\sum_{i=1}^k (b_i-a_i+1)}$ intervals of
length $2^{-a_{k+1}}$ each, which yields, by Proposition \ref{Proposition 4.1 falconer},
\begin{equation} \label{upper-bound for X'}
\begin{split}
 \dim_{\HH} X'(\a,\b)&\leq {\underline{\dim}}_{\B} X'(\a,\b)\leq 
    \liminf_{k\to\infty} \frac{\log 2^{\sum_{i=1}^k (b_i-a_i+1)}}
                              {-\log 2^{-a_{k+1}}}\\
      &=\liminf_{k\to\infty} \frac{\sum_{i=1}^k
(b_i-a_i+1)}{a_{k+1}}=d_1'.   
\end{split}
\end{equation}
For the lower bound, we now construct a measure which is somewhat more cumbersome. In fact, the measure $\mu$, constructed in the proof of Theorem \ref{hausdorff-dim-calc'}.1, is a convolution of measures, each of which ``corresponds'' to a single digit in the binary expansion of numbers in $[0,1]$. We now need to account for blocks of digits. With the unconstrained digits we deal as before (although now it is more convenient to group them in blocks, as we do for the constrained ones). However, the constrained
digits come in blocks, and we must regard them as basic units. Here we take the measure
$\mu'=*_{n=1}^\infty (\mu_n * \mu'_n)$, where
\begin{equation*}
 \begin{array}{rcl}
     \mu_n & =&\displaystyle\frac{1}{2^{b_n-a_n}} \sum_{l=0}^{2^{b_n-a_n}-1}
\delta_{l/2^{b_n}}\\
      &=& \displaystyle\left(\frac{1}{2}\delta_0 + \frac{1}{2}\delta_{2^{-(a_n +1)}}\right)*  \left(\frac{1}{2}\delta_0 + \frac{1}{2}\delta_{2^{-(a_n +2)}}\right)*\cdots*  \left(\frac{1}{2}\delta_0 + \frac{1}{2}\delta_{2^{-b_n}}\right),\\
      &&\\
    \mu'_n &= &\displaystyle \frac{1}{2} \delta_0 + \frac{1}{2}
\delta_{2^{-b_n}-2^{-a_{n+1}}}. 
 \end{array}   
\end{equation*}
Clearly, $\mu'$ is supported on $X'(\a,\b)$. Let $I= \left[\frac{l}{2^n}, \frac{l+1}{2^n}\right]$ for some $n$ and $0\leq l\leq 2^n-1$. We have
\begin{itemize}
\item If $l/2^n \notin \bigcup_{i=1}^\infty D_i'$, then
\begin{equation}\label{proof sec, equations1 for measure I for X'}
 \mu'(I)=0.   
\end{equation}
\item If $l/2^n \in  \bigcup_{i=1}^\infty D_i'$ and $a_{k}+ 1\leq n\leq b_{k}$ for some $k$, then
\begin{equation}\label{proof sec, equations3 for measure I for X'}
\mu'(I)= \displaystyle \frac{1}{2^{\sum_{i=1}^{k-1}(b_i-a_i+1) + (n-a_{k})}}.   \end{equation}
\item If $l/2^n \in \bigcup_{i=1}^\infty D_i'$ and $b_{k}+1\leq n\leq a_{k+1}$ for some $k$, then
\begin{equation}\label{proof sec, equations2 for measure I for X'}
\mu'(I)=  \displaystyle \frac{1}{2^{\sum_{i=1}^{k}(b_i-a_i+1)}}.   
\end{equation}
\end{itemize}
Let $\varepsilon>0$. The conditions of the theorem ensure that
\begin{equation}\label{proof sec, equation1 for inequality of sum in terms of eps X'}
\begin{split}   
\sum_{i=1}^{k-1} (b_i-a_i+1)\geq (d_1'-\varepsilon) a_{k},
\end{split}    
\end{equation}
for all sufficiently large $k$.
Hence, for $n\geq a_{k}+1$,
\begin{equation}\label{proof sec, equation2 for inequality of sum in terms of eps X'}
\begin{split}   
\sum_{i=1}^{k-1} (b_i-a_i+1) + (n-a_{k})&\geq n- (1- d_1'+\varepsilon)a_{k}\\
&\geq n- (1- d_1'+\varepsilon)n= (d_1'-\varepsilon)n.
\end{split}    
\end{equation}
By \eqref{proof sec, equations1 for measure I for X'}-\eqref{proof sec, equation2 for inequality of sum in terms of eps X'}, we have
\begin{equation}\label{proof sec, equation of measure of I' for X'}
    \mu'(I)\leq \displaystyle \frac{1}{2^{\sum_{i=1}^{k-1}(b_i-a_i+1) + (n-a_{k})}}\leq \frac{1}{2^{n(d_1'-\varepsilon)}}= |I|^{d_1'-\varepsilon}.
\end{equation}
Using \eqref{proof sec, equation of measure of I' for X'}, similarly to the proof of Theorem \ref{hausdorff-dim-calc'}, we obtain for every sufficiently small $I\subseteq [0,1]$,
\begin{equation}\label{proof sec, equation of measure of any I' for X'}
\mu'(I)\leq 2^{1+d_1'} |I|^{d_1'-\varepsilon}.    
\end{equation}
By \eqref{proof sec, equation of measure of any I' for X'} and Proposition \ref{Proposition mass distribution falconer},
$$\dim_{\HH} X'(\a,\b)\ge d_1'-\varepsilon,$$
and therefore
\begin{equation} \label{lower-bound for X'}
\dim_{\HH} X'(\a,\b)\ge d_1'.
\end{equation}
Using \eqref{pre, relation in box and modified box}, \eqref{pre, relation in Hausdoff and modified},  \eqref{upper-bound for X'}, and \eqref{lower-bound for X'},
we get
$$d_1'\leq \dim_{\HH} X'(\a,\b) \leq {\underline{\dim}}_{\MB} X'(\a,\b)\leq {\underline{\dim}}_{\B}X'(\a,\b) \leq d_1',$$
which completes the proof of this part.
\vspace{1.0 cm}

\noindent 2.
Denote 
$$d_2'=\limsup_{k\to\infty} \frac{\sum_{i=1}^k (b_i-a_i+1)}{b_{k}}.$$
We start with the rightmost equality in \eqref{proposition on X' equation of packing dim}.
For each $n$, cover $X'(\a,\b)$ by the union of all dyadic intervals $\left[\frac{l}{2^n}, \frac{l+1}{2^n}\right]$, $0\leq l\leq 2^n-1$, intersecting $X'(\a,\b)$. These intervals correspond to numbers $l$ whose binary expansion has, for each $i\geq 1$, the same digit in all places belonging to $[b_{i}+1, a_{i+1}]$. The number of these intervals is
\begin{equation}\label{equation upper box dimension Nn''}
\begin{split}
M_{2^{-n}}&= \begin{dcases}
   2^{(b_1-a_1 +1)+ \cdots + (b_{k-1}-a_{k-1} +1)+ (n- a_k)}, &\qquad a_{k}+ 1\leq n\leq b_k,\\
2^{(b_1-a_1 +1)+ \cdots + (b_k- a_k +1)}, &\qquad b_{k}+ 1\leq n\leq a_{k+1}.\\
\end{dcases}
\end{split}   
\end{equation}
Therefore, 
\begin{equation}\label{equation upper box dimension Nn/en1''}
\begin{split}
\frac{\log M_{2^{-n}}}{-\log 2^{-n}}
&= \begin{dcases}
    \frac{(b_1-a_1 +1)+ \cdots + (b_{k-1}-a_{k-1} +1)+ (n- a_k)}{n},  \quad &a_{k}+ 1\leq n\leq b_k,\\
  \frac{(b_1-a_1 +1)+ \cdots + (b_k- a_k +1)}{n}, \quad&b_{k}+ 1\leq n\leq a_{k+1}.\\
\end{dcases}\\
\end{split}   
\end{equation}
It follows from \eqref{equation upper box dimension Nn/en1''} that the sequence $\left(\frac{\log M_{2^{-n}}}{-\log 2^{-n}}\right)_{n=1}^\infty$ is increasing in the interval $[a_{k}+1, b_{k}+1]$ and decreasing in $[b_{k}+1, a_{k+1}]$. By Proposition~\ref{Proposition 4.1 falconer} and \eqref{equation upper box dimension Nn/en1''}, we have
\begin{equation}\label{equation upper box dimension bound from above''}
\begin{split}
{\overline{\dim}}_{\B} F &\leq \limsup_{n\to \infty} \frac{\log M_{2^{-n}}}{-\log 2^{-n}} = \limsup_{k\to \infty} \frac{\sum_{i=1}^k(b_i-a_i+1)}{b_k+1}\\
&= \limsup_{k\to \infty} \frac{\sum_{i=1}^k(b_i-a_i+1)}{b_k}=d_2'.    
\end{split}   
\end{equation}

In the other direction, consider the sets
\begin{equation*}
 \begin{split}
D_k''= \left\{x= \sum_{n=1}^{b_{k}}\frac{\xi_n}{2^n}\, :\, \xi_{b_i+1}=\xi_{b_i+2}=\cdots= \xi_{a_{i+1}},\, 1\leq i\leq k-1\right\},\qquad k\geq 1.  
 \end{split}   
\end{equation*}
We note that every two elements of $D_k''$ are at least $2^{-b_k}$ apart. Hence, each subinterval of $[0,1]$ of length $2^{-b_k}$ contains at most two points of $D_k''$. Since
$$|D_k''|=2^{\sum_{i=1}^{k-1}(b_i-a_i+1) + (b_k-a_k)},$$
the minimal number of intervals of length $2^{-b_k}$ required to cover $D_k''$ is at least 
$$\frac{1}{2}\cdot2^{\sum_{i=1}^{k-1}(b_i-a_i+1) + (b_k-a_k)}.$$
Let $N_{2^{-b_k}}$ be the minimal number of intervals of length $2^{-b_k}$ covering the set $X'(\a,\b)$. Since $D_k''\subseteq X'(\a,\b)$, we have 
\begin{equation}\label{equation for N(bk)''}
N_{2^{-b_k}}\geq \frac{1}{2}\cdot 2^{\sum_{i=1}^{k-1}(b_i-a_i+1) +(b_k-a_k)}.    
\end{equation}
By \eqref{equation for N(bk)''},
$$\frac{\log N_{2^{-b_k}}}{-\log 2^{-b_k}}\geq \frac{\log (2^{-1} 2^{\sum_{i=1}^{k-1}(b_i-a_i+1) +(b_k-a_k)})}{-\log 2^{-b_k}}.$$
Hence,
\begin{equation}\label{equation upper box dimension bound from below of Nn/en''}
\begin{split}
 {\overline{\dim}}_{\B} X'(\a,\b)&=\limsup_{\delta \to 0} \frac{\log N_{\delta}}{-\log \delta }\geq \limsup_{k\to \infty} \frac{\log N_{2^{-b_k}}}{-\log 2^{-b_k}}\\
 &\geq \limsup_{k\to \infty} \frac{\sum_{i=1}^{k-1}(b_i-a_i +1)+ (b_k-a_k)}{b_k}\\
 &= \limsup_{k\to \infty} \frac{\sum_{i=1}^{k}(b_i-a_i +1)}{b_k}=d_2'.
\end{split}
\end{equation}
By \eqref{equation upper box dimension bound from above''} and
\eqref{equation upper box dimension bound from below of Nn/en''} 
$$ {\overline{\dim}}_{\B} X'(\a,\b)=d_2'.$$
The two other equalities in \eqref{proposition on X' equation of packing dim} are proved similarly to the analogous equalities in \eqref{result sec, formula for packing dim of X i}.
\end{Proof}

\section{Proof of Theorem \ref{E_0-may-be-any-dim}}\label{section proof A}
\medskip
Let $\a=(a_i)_{i=1}^\infty$ and $\b=(b_i)_{i=1}^\infty$ be as
in Theorem \ref{hausdorff-dim-calc'}, and satisfying,
moreover, the conditions
\begin{equation}
a_i+i<b_i, \qquad i=1,2,\ldots,
\end{equation}
and
\begin{equation}
\frac{b_i}{a_i} \xrightarrow[i\to\infty]{} \infty,
   \qquad \frac{b_i}{a_{i+1}} \xrightarrow[i\to\infty]{} d\,.
\end{equation}
Let $\r$ be the sequence consisting of all
numbers of the form $2^n$, where $b_{i-1}+1\leq n\le a_i$ for some $i$.

\noindent 1. We divide the proof into two parts:
\begin{description}
\item {$(i)$} \vspace{-0.3 cm} $$\dim_{\HH} E_{0}(\r)= d.$$
\item {$(ii)$}  \vspace{-0.3 cm} $${\underline{\dim}}_{\MB} E_{0}(\r)= d.$$
\end{description}
\vspace{0.25 cm}

\noindent $(i)$ We first claim that 
\begin{equation}\label{proof sec, equation of X(a,b) subset of E0(r)}
E_0(\r)\supseteq X(\a',\b),    
\end{equation}
where $\a'=(a_i+i)_{i=1}^\infty$.
In fact, take a typical element
of $\r$, say $2^{n_0}$,~where $b_{i_0-1}+1\leq n_0\le a_{i_0}$ for some
$i_0$, and $x\in X(\a',\b)$. Let $x=\sum_{n=1}^\infty\xi_n
2^{-n}$. Then:
\begin{equation} \label{separation}
\begin{split}
2^{n_0} x &=\sum_{n=1}^{n_0} \xi_n 2^{n_0-n}
             +2^{n_0}\sum_{n=n_0+1}^{a_{i_0}+ i_0} \xi_n 2^{-n}
             +2^{n_0}\sum_{n=a_{i_0}+ i_0 +1}^\infty \xi_n 2^{-n}.\\ 
\end{split}
\end{equation}
The first term on the right-hand side is an integer, the second consists of zeros, and the
third is bounded above by $2^{-i_0}$. Hence $2^n x\longto 0\;
(\textup{mod}\, 1)$ as $n\to\infty$ with $2^n$ belonging to $\r$, and therefore $E_0(\r)\supseteq X(\a',\b)$. Using Theorem \ref{hausdorff-dim-calc'}.1 and the monotonicity of the Hausdorff dimension, we get
\begin{equation}\label{proof sec, equation of Dim(H)(E0)>d}
\begin{split}
 \dim_{\HH} E_0(\r)&\ge \dim_{\HH} X(\a',\b)=  {\underline{\dim}}_{\B} X(\a',\b) =  \liminf_{k\to\infty} \frac{\sum_{i=1}^k
(b_i-a_i')}{a_{k+1}'}\\
&= \liminf_{k\to\infty} \frac{\sum_{i=1}^k
(b_i-a_i-i)}{a_{k+1}+(k+1)}=\lim_{k\to \infty}\frac{b_{k}}{a_{k+1}}=d.
\end{split}    
\end{equation}

On the other hand, if $x=\sum_{n=1}^\infty \xi_n 2^{-n} \in E_0(\r)$, then for sufficiently large $n$ we have in particular $r_n x\in (-1/4,1/4)\; (\textup{mod}\, 1)$. Now, if $2^l x\in (-1/4,1/4) \; (\textup{mod}\, 1)$ for some $l$, then
$\xi_{l+1}=\xi_{l+2}$. Since, for sufficiently large $i$, all
numbers
$$2^{b_i+1} x, 2^{b_i+2} x,\ldots,2^{a_{i+1}} x$$
lie in
the interval $(-1/4,1/4)$ modulo~1, it means that all digits
$$\xi_{b_i+2},\xi_{b_i+3},\ldots,\xi_{a_{i+1}+2}$$
are equal.
Thus,
\begin{equation}\label{proof sec, equation E0(r) subset of union of X'(a,b)}
E_0(\r)\subseteq \bigcup_{j=1}^\infty X'(\a''(j),\b'(j)),    
\end{equation}
where
$\a''(j)=(a_i+2)_{i=j}^\infty$ and $\b'(j)=(b_i+1)_{i=j}^\infty$.
Note that
\begin{equation}\label{proposition on X' equation on X' subset of X'}
 X'(\a''(1),\b'(1))\subseteq X'(\a''(2),\b'(2))\subseteq \cdots.   
\end{equation}
By \eqref{Hausdorff dim of union}, Theorems \ref{hausdorff-dim-calc''}.1, \eqref{proof sec, equation E0(r) subset of union of X'(a,b)}, and \eqref{proposition on X' equation on X' subset of X'},
\begin{equation}\label{proof sec, equation of Dim(H)(E0)<d}
\begin{split}
\dim_{\HH} E_0(\r)&\leq \dim_{\HH} \bigcup_{j=1}^\infty X'(\a''(j),\b'(j))=\sup\left\{\dim_{\HH} X'(\a''(j),\b'(j))\, :\, j=1,2,\ldots\right\}\\
&=\lim_{j\to \infty}\dim_{\HH} X'(\a''(j),\b'(j))=\lim_{j\to \infty}\left\{ \liminf_{k\to\infty} \frac{\sum_{i=j}^k (b_i'-a_i''+1)}{a_{k+1}'}\right\}\\
&=\lim_{j\to \infty}\left\{ \liminf_{k\to\infty} \frac{\sum_{i=j}^k (b_i-a_i)}{a_{k+1}+2}\right\}\leq\liminf_{k\to\infty} \frac{\sum_{i=1}^k (b_i-a_i)}{a_{k+1}+2}\\
&\leq\liminf_{k\to\infty} \frac{\sum_{i=1}^k (b_i-a_i)}{a_{k+1}}= d.
\end{split}   
\end{equation}
By \eqref{proof sec, equation of Dim(H)(E0)>d} and \eqref{proof sec, equation of Dim(H)(E0)<d}
$$\dim_{\HH} E_0(\r)=d,$$
which completes the proof in this part.
\vspace{1 cm}

\noindent $(ii)$ Using \eqref{proof sec, equation of X(a,b) subset of E0(r)} and \eqref{proof sec, equation E0(r) subset of union of X'(a,b)}, we have 
\begin{equation}\label{proof sec, E0 in between two sets}
X(\a',\b)\subseteq E_0(\r)\subseteq \bigcup_{j=1}^\infty X'(\a''(j),\b'(j)).    
\end{equation}
Clearly, $X(\a',\b)\subseteq X'(\a''(1),\b'(1))$. Using \eqref{proposition on X' equation on X' subset of X'}, we conclude that 
\begin{equation}\label{proof sec, X is containd in X'}
X(\a',\b)\subseteq X'(\a''(j), \b'(j)),\qquad j\geq 1.
\end{equation}
It follows from Theorems \ref{hausdorff-dim-calc'}-\ref{hausdorff-dim-calc''}, \eqref{proof sec, equation of Dim(H)(E0)<d}, and  \eqref{proof sec, X is containd in X'} that
\begin{equation}\label{proof sec, equation of Hausdorff and lower box for X and X' are same}
\begin{split}
\dim_{\HH} X(\a',\b)&={\underline{\dim}}_{\MB} X(\a',\b)= {\underline{\dim}}_{\B} X(\a',\b)= \dim_{\HH} X'(\a''(j), \b'(j))\\
&= {\underline{\dim}}_{\MB} X'(\a''(j), \b'(j))={\underline{\dim}}_{\B} X'(\a''(j), \b'(j))=d,    \end{split}\qquad j\geq 1.
\end{equation}
By \eqref{proof sec, E0 in between two sets} and \eqref{proof sec, X is containd in X'}, we have
\begin{equation}\label{proof sec, E0 equal to infinte union}
X(\a',\b)\subseteq E_0(\r)=\bigcup_{j=1}^\infty \left(X'(\a''(j), \b'(j))\bigcap E_0(\r)\right),    
\end{equation}
and 
\begin{equation}\label{proof sec, X is containd in X' for each j}
X(\a',\b)\subseteq X'(\a''(j),\b'(j))\bigcap E_0(\r),\qquad j\geq 1.
\end{equation}
Using the monotonicity of the lower box dimension, \eqref{proof sec, equation of Hausdorff and lower box for X and X' are same}, and  \eqref{proof sec, X is containd in X' for each j}, we obtain
$$d={\underline{\dim}}_{\B}X(\a',\b)\leq {\underline{\dim}}_{\B} X'(\a''(j), \b'(j))\bigcap E_0(\r)\leq {\underline{\dim}}_{\B}X'(\a''(j), \b'(j)) =d, \quad j\geq 1.$$
Hence 
\begin{equation}\label{proof sec, box dim of X' intersect with E_0}
{\underline{\dim}}_{\B}X'(\a''(j),\b'(j))\bigcap E_0(\r)=d,\qquad j\geq 1.  
\end{equation}
Since $\dim_{\HH} X(\a',\b)=d$, by \eqref{pre, relation in Hausdoff and modified}, we have
\begin{equation}\label{proof sec, modified lower box of X(a;,b) bounded below}
d= \dim_{\HH} X(\a',\b)\leq {\underline{\dim}}_{\MB}X(\a',\b).    
\end{equation}
Using the monotonicity of the modified lower box dimension, \eqref{proof sec, E0 equal to infinte union}, and \eqref{proof sec, modified lower box of X(a;,b) bounded below}, we get 
\begin{equation}\label{proof sex, modified lower box of E_0 bounded below}
\begin{split}
d\leq {\underline{\dim}}_{\MB}X(\a',\b) \leq  {\underline{\dim}}_{\MB}E_{0}(\r).   
\end{split} 
\end{equation}
By \eqref{pre, relation in box and modified box}, the countable stability of the  modified lower box dimension, \eqref{proof sec, E0 equal to infinte union}, \eqref{proof sec, box dim of X' intersect with E_0}, and \eqref{proof sex, modified lower box of E_0 bounded below}, we have
\begin{equation*}
\begin{split}
d\leq {\underline{\dim}}_{\MB}E_0(\r)&={\underline{\dim}}_{\MB} \bigcup_{j=1}^\infty \left(X'(\a''(j), \b'(j))\bigcap E_0(\r)\right)\\
&= \sup\left\{{\underline{\dim}}_{\MB} \left(X'(\a''(j), \b'(j))\bigcap E_0(\r)\right)\, :\, j=1,2,\ldots\right\}\leq d.
\end{split}   
\end{equation*}
Hence,
$${\underline{\dim}}_{\MB}E_0(\r)=d.$$
\vspace{1 cm}

\noindent 2. We divide the proof into two parts:
\begin{description}
\item {$(i)$} \vspace{-0.3 cm} $$\dim_{\HH} E_{\f}(\r)= d.$$
\item {$(ii)$}  \vspace{-0.3 cm} $${\underline{\dim}}_{\MB} E_{\f}(\r)= d.$$
\end{description}
\vspace{0.25 cm}

\noindent $(i)$ Let $\a=(a_i)_{i=1}^\infty,\b=(b_i)_{i=1}^\infty$, $\a'=(a_i+i)_{i=1}^\infty$ and
$\r=(r_n)_{n=1}^\infty$ be as
in the proof of part~$(1)$ of the theorem.
Let $x\in X(\a',\b)$. Similarly to (\ref{separation}), note
that
\begin{equation*}
\begin{split}
\sum_{n=1}^\infty \|r_n x\|
&=\sum_{k=1}^\infty \sum_{n=b_k+1}^{a_{k+1}} \|2^n x\|
\leq \sum_{k=1}^\infty \sum_{n=b_k+1}^{a_{k+1}} \frac{1}{2^{a_{k+1}'-n}}=\sum_{k=1}^\infty \sum_{n=b_k+1}^{a_{k+1}} \frac{1}{2^{a_{k+1}+ k+1-n}}\\
&= \sum_{k=1}^\infty \sum_{j=0}^{a_{k+1}-b_k-1} \frac{1}{2^{k+1+j}}< \sum_{k=1}^\infty \frac{1}{2^{k}}= 1,
\end{split}
\end{equation*}
so that $E_\f(\r)\supseteq X(\a',\b)$, and therefore
$\dim_{\HH} E_\f(\r)\ge d$.

In the other direction, since $E_\f(\r)\subseteq E_0(\r)$ and $\dim_{\HH} E_0(\r)=d$, we have
$$\dim_{\HH} E_\f(\r)\le \dim_{\HH} E_0(\r)= d.$$
Hence,
$$\dim_{\HH} E_\f(\r)=d.$$
\vspace{0.25 cm}

\noindent $(ii)$ As above, we have
\begin{equation}\label{proof sec, eqution X in Ef and Ef in E0}
 X(\a',\b)\subseteq E_\f(\r)\subseteq E_0(\r),   
\end{equation}
and 
\begin{equation}\label{proof sec, hausdorff dim of X(a',b) equal to modifed of E_0}
\dim_{\HH} X(\a',\b)=d,\qquad {\underline{\dim}}_{\MB}E_0(\r)=d.    
\end{equation}
By \eqref{pre, relation in Hausdoff and modified}, \eqref{proof sec, eqution X in Ef and Ef in E0}, \eqref{proof sec, hausdorff dim of X(a',b) equal to modifed of E_0}, and the monotonicity of the modified lower box dimension, we get
$$d= \dim_{\HH} X(\a',\b)\leq {\underline{\dim}}_{\MB}X(\a',\b)\leq {\underline{\dim}}_{\MB}E_\f(\r)\leq {\underline{\dim}}_{\MB}E_0(\r)=d.$$
Hence,
$${\underline{\dim}}_{\MB}E_\f(\r)=d,$$
which completes the proof of this part.
\vspace{1 cm}

\section{IP-sets; Proof of Theorem \ref{E-may-be-any-dim}}\label{Section on IP sets}
\medskip
Before getting to the proof of Theorem \ref{E-may-be-any-dim}, we recall a few relevant results regarding IP-sets. A set~$P$ of positive integers is an {\it
IP-set} (cf.\ \cite{furstenberg}) if, for some sequence
$\p=(p_n)_{n=1}^\infty$ of positive integers,~$P$ consists of all
finite sums of the form $p_{i_1} + p_{i_2} + \cdots + p_{i_k}$ with
$i_1 < i_2 < \cdots < i_k$.  The elements of $P$ can be ordered in a natural way to form a sequence as follows. If $l=\sum_{k=1}^d \varepsilon_k 2^{k-1}$ is the binary representation of $l$, then the $l$-th term of our sequence is $\sum_{k=1}^d\varepsilon_k p_k$. (Note that the sequence may assume some values more than once.) We will denote this set by IP-$(p_n)$ and refer to it as the {\it IP-set generated by} $(p_n)$.
 
Similarly, we can define IP-sets in any additive semigroup. Here, we will be interested also in IP-sets in the circle group $\T=\R/\Z$. Let $(\alpha_n)_{n=1}^\infty$ be a sequence in $\T$. We recall \cite[Proposition 2.1]{berend} that a necessary condition for IP-$(\alpha_n)$ to be dense in $\T$ is
\begin{equation}
\sum_{n=1}^\infty \|h\alpha_n\|\ge 1, \qquad h=1,2,\ldots,
\end{equation}
while the stronger condition
\begin{equation}\label{proof sec, equation alpha not dense}
\sum_{n=1}^\infty \|h\alpha_n\|=\infty, \qquad h=1,2,\ldots,
\end{equation}
is sufficient for the same to hold.
\vspace{1 cm}

\begin{Proof}{\ of Theorem \ref{E-may-be-any-dim}} 
Let $\a=(a_i)_{i=1}^\infty$, $\b=(b_i)_{i=1}^\infty$, and  $\a'=(a_i+i)_{i=1}^{\infty}$ be as in the proof of Theorem \ref{E_0-may-be-any-dim}. Let $\p$ be the subsequence  of $(2^n)_{n=1}^\infty$, consisting of all powers $2^n$ with $n\in \bigcup_{i=1}^\infty [b_{i}+1, a_{i+1}]$. Let $(r_n)_{n=1}^\infty$ be the IP-sequence generated by $\p$. For $x\in \R$, denote by $\mathcal{S}_{x}$ the IP-set generated in $\T$ by $\{2^n x\; (\textup{mod}\, 1)\, :\, n\in \bigcup_{i=1}^\infty [b_{i}+1, a_{i+1}]\}$. Note that $\mathcal{S}_{x}=\{r_n x\; \; (\textup{mod}\, 1) \, :\, n\geq 1\}$.

We divide the proof into two parts:
\begin{description}
\item {$(i)$} \vspace{-0.3 cm} $$\dim_{\HH} E(\r)= d.$$
\item {$(ii)$}  \vspace{-0.3 cm} $${\underline{\dim}}_{\MB} E(\r)= d.$$
\end{description}
\vspace{0.25 cm}

\noindent $(i)$ We first claim that $X(\a',\b)\subseteq E(\r)$.    
In fact, let $x=\sum_{n=1}^\infty \xi_n 2^{-n}\in~X(\a',\b)$. We have to show that $\mathcal{S}_{x}$ is not dense in $\T$. Let $b_{i_0-1}+1\leq  n_0\leq a_{i_0}$ for some~$i_0$. We have 
\begin{equation}\label{proof sec, equation for 2^n0}
\begin{split}
2^{n_0}x &= \sum_{n=1}^{n_{0}}\xi_n 2^{n_0-n} + 2^{n_0}\sum_{n= n_0+1}^{a_{i_0}'} \xi_n 2^{-n}+ 2^{n_0}\sum_{n= a_{i_0}'+1}^\infty \xi_n 2^{-n}.\\
\end{split}    
\end{equation}
The first sum on the right-hand side of \eqref{proof sec, equation for 2^n0} is an integer, and the second consists of zeros only. Hence
$$2^{n_0}x\equiv 2^{n_0}\sum_{n= a_{i_0}'+1}^\infty \xi_n 2^{-n} \; (\textup{mod}\, 1).$$
It follows that:
$$\|2^{n_0}x\| \leq \frac{1}{2^{a_{i_0}+i_0-n_0}}.$$
\end{Proof}
We have 
\begin{equation}
\begin{split}
\sum_{n\in \bigcup_{i=1}^\infty[b_{i}+1,a_{i+1}]}\|2^n x\|&\leq \sum_{n\in \bigcup_{i=1}^\infty[b_{i}+1,a_{i+1}]}\frac{1}{2^{a_{i+1}+ i+1-n}}\\
&=\sum_{i=1}^\infty \sum_{j=0}^{a_{i+1}-b_i-1}\frac{1}{2^{i+1+j}}\\
&< \sum_{i=1}^\infty \frac{1}{2^{i+1}}\sum_{j=0}^\infty\frac{1}{2^{j}}=\sum_{i=1}^\infty \frac{1}{2^{i}}=1,
\end{split}    
\end{equation}
which implies that $\mathcal{S}_{x}$ is not dense modulo 1. Hence 
\begin{equation}\label{proof sec, X(a,b) contained in E(r)}
X(\a',\b)\subseteq E(\r).    
\end{equation}
By \eqref{proof sec, equation of Hausdorff and lower box for X and X' are same},
\begin{equation}\label{proof sec, dim E(r)>d}
\dim_{\HH}E(\r)\geq \dim_{\HH}X(\a',\b)=d.    
\end{equation}

In the other direction, we will show that $E(\r)$ is contained in some set of Hausdorff dimension $d$. Let $x= \sum_{n=1}^\infty \xi_n2^{-n} \in E(\r)$. Then $\mathcal{S}_{x}$ is not dense modulo 1. By~\eqref{proof sec, equation alpha not dense}, there exists an $h\in \N$ such that 
$$\sum_{n\in \bigcup_{i=1}^\infty [b_{i}+1, a_{i+1}]}\|h 2^n x\|<\infty.$$
which implies that, as $n\to \infty$ along the set $\bigcup_{i=1}^\infty [b_{i}+1, a_{i+1}]$,
$$\lim\|h2^n x\|=0.$$
We have $h2^n x\in \left(-\frac{1}{4},\frac{1}{4}\right)\, (\textup{mod}\, 1)$ for all sufficiently large $n$. Using a similar calculation to that in the proof of Theorem \ref{E_0-may-be-any-dim}, we get $hx\in X'(\a''(i),\b'(i))$, where $\a''(i)= (a_j+2)_{j=i}^\infty$ and $\b'(i)=(b_j+1)_{j=i}^\infty$ for some $i\geq 1$.
Hence
\begin{equation}\label{proof sec, E(r) contained in union of X'(a,b)}
E(\r)\subseteq \bigcup_{h=1}^\infty \bigcup_{i=1}^\infty h^{-1} X'(\a''(i), \b'(i)).    
\end{equation}
Here, for $A\subseteq \T$, we have denoted by $h^{-1}A$ the set $\{x\in \T\, :\, hx \in A\}$. Note that, by the stability and Lipschitz invariance of Hausdorff dimension, 
\begin{equation}\label{proof sec, equation of dim(H)(1/h)A= dim(H)A}
\dim_{\HH}h^{-1}A = \dim_{\HH}A   
\end{equation}
for every set $A\subseteq \T$.
Denote 
\begin{equation}\label{proof sec, equation of A[h] expression}
A[h]= h^{-1}\bigcup_{i=1}^\infty X'(\a''(i), \b'(i)),\qquad h\geq 1. \end{equation}
It follows from the proof of Theorem \ref{E_0-may-be-any-dim} that
$\dim_{\HH} A[1]= d$. By \eqref{proof sec, equation of dim(H)(1/h)A= dim(H)A}, we have
\begin{equation}\label{proof sec, equation of Hausdorff dim of A[h]}
  \dim_{\HH}A[h]= d,\qquad h\geq 1.  
\end{equation}
By \eqref{proof sec, E(r) contained in union of X'(a,b)}, \eqref{proof sec, equation of A[h] expression}, and \eqref{proof sec, equation of Hausdorff dim of A[h]},

\begin{equation}\label{proof sec, dim E(r)<d}
\begin{split}
\dim_{\HH} E(\r)&\leq \dim_{\HH}\bigcup_{h=1}^\infty \bigcup_{i=1}^\infty h^{-1} X'(\a''(i), \b'(i)) = \dim_{\HH}\bigcup_{h=1}^\infty A[h]\\  
&=\sup\{\dim_{\HH}A[h]\, :\, h=1,2,\ldots\}= d.
\end{split}  
\end{equation}
By \eqref{proof sec, dim E(r)>d} and \eqref{proof sec, dim E(r)<d}, we have
$$\dim_{\HH} E(\r)=d,$$
which completes the proof of this part.
\vspace{1 cm}

\noindent $(ii)$ By  \eqref{proof sec, X(a,b) contained in E(r)}, \eqref{proof sec, E(r) contained in union of X'(a,b)}, and \eqref{proof sec, equation of A[h] expression}, we have
\begin{equation}\label{proof sec, E(r) contained in two sets}
X(\a',\b)\subseteq E(\r)\subseteq \bigcup_{h=1}^\infty \bigcup_{i=1}^\infty h^{-1} X'(\a''(i), \b'(i))= \bigcup_{h=1}^\infty A[h].
\end{equation}
By \eqref{proof sec, equation of Hausdorff and lower box for X and X' are same},
\begin{equation}\label{proof sec, equation modified box of X(a',b)}
 {\underline{\dim_{\B}}} X(\a',\b)={\underline{\dim_{\MB}}} X(\a',\b)=d,   
\end{equation}
and 
\begin{equation}\label{proof sec, lower box for X(a'',b')<d}
   {\underline{\dim_{\MB}}} X'(\a''(i), \b'(i))=d,\qquad i\geq 1.
\end{equation}
By the stability and Lipschitz invariance of the modified lower box dimension, we have
\begin{equation}\label{proof sec, equation of dim(MB)(1/h)A= dim(MB)A}
 {\underline{\dim_{\MB}}}h^{-1}A = {\underline{\dim_{\MB}}}A   
\end{equation}
for every set $A\subseteq \T$.
Hence,
\begin{equation}\label{proof sec, equatiom of modified of A[h] equal to A[1]}
{\underline{\dim_{\MB}}} A[h]= {\underline{\dim_{\MB}}} A[1],\qquad h\geq 1.    
\end{equation}
By \eqref{proof sec, E(r) contained in two sets}, \eqref{proof sec, equation modified box of X(a',b)}, \eqref{proof sec, lower box for X(a'',b')<d}, and \eqref{proof sec, equatiom of modified of A[h] equal to A[1]},
\begin{equation}
\begin{split}
d&= {\underline{\dim}}_{\MB}X(\a',\b)
\leq {\underline{\dim}}_{\MB} E(\r)\leq {\underline{\dim}}_{\MB}\bigcup_{h=1}^\infty A[h]\\
&= \sup\left\{{\underline{\dim}}_{\MB}A[h]\,:\, h=1,2,\ldots\right\}={\underline{\dim}}_{\MB}A[1]\\
&= \sup\left\{ {\underline{\dim_{\MB}}} X'(\a''(i), \b'(i))\,:\, i=1,2,\ldots\right\}=d.
\end{split}   
\end{equation}
Hence,
$${\underline{\dim}}_{\MB} E(\r)=d.$$
\nopagebreak

\bigskip
\end{document}